# ESTIMATES FOR NORMS OF RANDOM POLYNOMIALS

PAVEL G. GRIGORIEV

*A-1611, Dept. of Mechanics and Mathematics,*
*Moscow State University, Moscow, 119899, RUSSIA*
E-mail: pgrigoriev@yahoo.com

ABSTRACT. This paper contains some estimates for the *integral-uniform* norm and the uniform norm of a wide class of random polynomials. The family of integral-uniform norms introduced in [6] is a natural generalization of the maximum norm taken over a net. We prove some properties of the integral-uniform norms. The given application of the established estimates demonstrates that the integral-uniform norms may be useful whenever one is interested in the properties of a function distribution.

## 1. INTRODUCTION

In this paper some estimates for mathematical expectation of norms of random polynomials of the type

$$\sum_{j=1}^{n} a_j \xi_j(\omega) f_j(x) \tag{1}$$

are presented. Here $\{\xi_i\}_1^n$ is a set of independent random values defined on $(\Omega, \mathsf{P})$ and $\{f_i\}_1^n$ is a set of functions on another probability space $(X, \mu)$. The norms here are taken in a space of functions, which depend only on the space variable $x$ with fixed $\omega$.

Similar estimates for various systems of functions $\{f_i\}_1^n$ and random variables $\{\xi_i\}_1^n$ have been widely applied in analysis since 1930s. In 1954 Salem and Zygmund [14] established a number of estimates for the uniform norm of random trigonometric polynomials. In particular in [14], it was shown that

$$\mathsf{E} \| \sum_{k=-n}^{n} r_k(\omega) e^{ikt} \|_\infty \asymp (n \log n)^{1/2},$$

where $r_k(\omega)$ are the Rademacher functions, here and further the expression $A_n \asymp B_n$ stands for $cA_n \leq B_n \leq CA_n$ with some constants $c$, $C$. This estimate along with Khinchin's inequality reflects subtle differences between a finite dimensional subspace of $L_\infty$ and its natural embeddings in $L_p$ spaces with $1 \leq p < \infty$.

By now various methods for estimating the uniform norm of random polynomials (1) have been developed (e.g. see [5], [8], [9]). A first lower estimate for the uniform norm

Almost the same text was published in East J. Approximation 2001 v. 7 no. 4 pp. 445–469.

This work was supported by the Grant *"Leading Scientific Schools of the Russian Federation"* 00-15-96047.

of a random polynomial (1) with respect to a *general function system* was established by Kashin and Tzafriri in [5]–[7], this result will be formulated in Section 2.

In [6] Kashin and Tzafriri introduced the following norm

$$\|f\|_{m,\infty} := \int_X \ldots \int_X \max\{|f(x_1)|, \ldots, |f(x_m)|\} d\mu(x_1)\ldots d\mu(x_m), \tag{2}$$

where $f$ is a function defined on a measure space $(X, \mu)$, $\mu(X) = 1$. This norm is a natural generalization of $\|\cdot\|_\infty$–norm taken over a net, we call it *the integral-uniform norm*. One can easily see that for every integrable function $f \in L_1(X)$ we get $\|f\|_1 = \|f\|_{1,\infty}$ and

$$\|f\|_{m,\infty} = \int_0^\infty \left(1 - (1 - \lambda_f(t))^m\right) dt, \tag{3}$$

where

$$\lambda_f(t) := \mu\{\tau : |f(\tau)| > t\}.$$

It is also easy to notice, that for $f \in L_\infty(X)$ the following inequalities take place $\|f\|_1 \leq \|f\|_{m,\infty} \leq \|f\|_\infty$ and $\|f\|_{m,\infty} \to \|f\|_\infty$ as $m \to \infty$. Using a trivial inequality $\max(|a|, |b|) \leq |a| + |b|$ and the definition of the integral-uniform norm (2), we get

$$\|f\|_{n,\infty} \leq \left(\frac{n}{m} + 1\right) \|f\|_{m,\infty}, \qquad m < n \tag{4}$$

for all $f \in L_1(X)$. For the integral-uniform norm of an indicator $\chi_\Delta$ of a set $\Delta \subset X$ the identity (3) implies

$$\|\chi_\Delta\|_{m,\infty} = 1 - (1 - \mu\Delta)^m, \tag{5}$$

Thus, if we take $m$ of order $1/\mu\Delta$ then $\|\chi_\Delta\|_{m,\infty}$ is of order one.

The technique used in [5], [7] for estimating the uniform norm of random polynomials turned out to be applicable for estimating the *integral-uniform* norm (2). In fact, an estimate for *the integral-uniform* norm of random polynomials (1) for a special case of parameter $m$ was implicitly obtained in [5]–[7].

In Section 2 we present some generalizations of the results from [5]–[7] for both the case of the $\|\cdot\|_{m,\infty}$–norm with an arbitrary parameter $m$ and a wider class of function systems $\{f_i\}_1^n$. The generalizations are obtained by the same method as in [5]–[7], which relies on a multidimensional version of the central limit theorem with precise estimate of the error term. In Section 3 we shall show that under some additional constraints on $\{\xi_i\}_1^n$ and $m$ the established estimate is precise in sense of order. In Section 4 we mention some properties of the integral-uniform norm, in particular, its properties are illustrated on some inequalities for the integral-uniform norms of trigonometric polynomials. In addition, in Section 4 we present an application of the established estimates. This application uses a simple geometrical lemma which could be of independent interest. Most of the results have presented here been announced by the author in [3].

I would like to express my special thanks to B.S. Kashin for his numerous useful comments and advices, also I am very grateful to E.M. Semenov for interesting discussions.

## 2. The lower estimates for the integral-uniform norms of random polynomials.

In [6], [7] Kashin and Tzafriri proved that whenever systems of functions $\{f_i\}_{i=1}^n$ and $\{\xi_i\}_{i=1}^n$, defined on probability spaces $(X, \mu)$ and $(\Omega, \mathsf{P})$ respectively, satisfy the following conditions

(a) $\|f_i\|_2 = 1$ and $\|f_i\|_3 \leq M$ for every $i$;
(b) $\|\sum_{i=1}^n c_i f_i\|_2 \leq M \left(\sum_{i=1}^n |c_i|^2\right)^{1/2}$ for all sets of coefficients $\{c_i\}_1^n$;
(c) $\{\xi_i\}_1^n$ are independent variables, such that $\mathsf{E}\xi_i = 0$, $\mathsf{E}|\xi_i|^2 = 1$ and $(\mathsf{E}|\xi_i|^3) \leq M^3$.

Then there exist positive constants $q = q(M)$, $C_j = C_j(M)$, $j = 1, 2, 3$ such that

$$(6) \quad \mathsf{P}\Big\{\omega \in \Omega : \|\sum_{i=1}^n a_i \xi_i(\omega) f_i\|_{L_\infty(X)} \leq C_1 \Big(\sum_{i=1}^n |a_i|^2\Big)^{\frac{1}{2}} (1 + \log R)^{\frac{1}{2}}\Big\} \leq \frac{C_2}{R^q},$$

where

$$(7) \quad R := \frac{(\sum_{i=1}^n |a_i|^2)^2}{\sum_{i=1}^n |a_i|^4},$$

and hence

$$(8) \quad \mathsf{E}\|\sum_{i=1}^n a_i \xi_i f_i\|_{L_\infty(X)} \geq C_3 \Big(\sum_{i=1}^n |a_i|^2\Big)^{1/2} (1 + \log R)^{1/2}.$$

The proof of these estimates practically involved the estimate of $\|\sum_1^n a_i \xi_i f_i\|_{m,\infty}$ for a special value of parameter $m$, precisely $m \asymp (1 + \log R)^2 R^{1/2+\varepsilon}$.

In this paper a generalization of the inequalities (6), (8) for both the case of integral-uniform norm and a wider class of random polynomials is established. In particular, it is shown that if $R(\{a_i\}_1^n) \asymp n$ then the estimates (6), (8) stay true whenever functions $\{f_i\}_1^n$ satisfy instead of (b) the following condition

(b') $\|\sum_{i=1}^n c_i f_i\|_2 \leq M n^p \left(\sum_{i=1}^n |c_i|^2\right)^{1/2}$ for all sets of coefficients $\{c_i\}_1^n$ with some constants $M > 0$ and $p \in [0, 1/2)$.

**Theorem 1.** *Let $\{f_i\}_{i=1}^n$ and $\{\xi_i\}_{i=1}^n$ be sets of functions defined on probability spaces $(X, \mu)$ and $(\Omega, \mathsf{P})$ respectively, which satisfy (a) and (c). Let also $\{a_i\}_1^n$ be a fixed set of coefficients and for all sets of coefficients $\{c_i\}_1^n$ the following inequality hold*

$$\|\sum_{i=1}^n c_i f_i\|_2 \leq M \big[R(\{a_i\})\big]^p \Big(\sum_{i=1}^n |c_i|^2\Big)^{1/2},$$

*where $R = R(\{a_i\})$ defined by (7); $M > 0$ and $p \in [0, 1/2)$ are some constants. Then there exist positive constants $q' = q'(p)$, $C'_j = C'_j(p, M)$, $j = 1, 2, 3$ such that the following estimates take place*

$$(9) \quad \mathsf{P}\Big\{\omega \in \Omega : \|\sum_{i=1}^n a_i \xi_i(\omega) f_i\|_{m,\infty} \leq C'_1 \Big(\sum_{i=1}^n |a_i|^2 \log P\Big)^{1/2}\Big\} \leq \frac{C'_2}{P^{q'}},$$

$$(10) \quad \mathsf{E}\|\sum_{i=1}^n a_i \xi_i f_i\|_{m,\infty} \geq C'_3 \Big(\sum_{i=1}^n |a_i|^2 \log P\Big)^{1/2},$$

where $P := \min(m, R) + 1$.

**Corollary 1.** *Let the coefficients $\{a_i\}_1^n$ satisfy $R(\{a_i\}_1^n) \asymp n$ and functions $\{f_i\}_1^n$ and $\{\xi_i\}_1^n$ satisfy* **(a)** *and* **(c)** *respectively. Let also $\{f_i\}_1^n$ satisfy*

$$\|\sum_{i=1}^n \varepsilon_i f_i\|_2 \leq M n^{\frac{1}{2}+p}$$

*for all signs $\varepsilon_i = \pm 1$ with some constants $p \in [0, 1/2)$, $M > 0$. Then the estimates (9), (10) hold for the random polynomial (1).*

To prove the Corollary it suffices to notice that Lemma from Sec. 4 implies the condition **(b′)** (with another $p \in [0, \frac{1}{2})$) for the functions $\{f_i\}_1^n$.

The proof of Theorem 1 essentially follows the pattern of the proof of (6), (8) from [7]. In Section 3 it is shown that under some additional constraints on $\xi_k$ and $m$ the estimate (10) is precise in sense of order.

**Proof of Theorem 1.**

**Step 1.** We can re-scale the coefficients $\{a_i\}_{i=1}^n$ so that $\sum_1^n |a_i|^2 = 1$. Let $\varepsilon = \varepsilon(M) = \frac{1}{4}(\frac{3}{4M^2})^3$. Consider the set

$$E_1 := \Big\{ x : \sum_{i=1}^n |a_i|^3 |f_i(x)|^3 < \frac{M^3}{\varepsilon} \sum_{i=1}^n |a_i|^3 \Big\}.$$

Then assumption **(a)** and Chebyshev's inequality imply

$$\mu(E_1^c) \frac{M^3}{\varepsilon} \sum_{i=1}^n |a_i|^3 \leq \int_\Omega \sum_{i=1}^n |a_i|^3 |f_i(x)|^3 d\mu(x) \leq M^3 \sum_{i=1}^n |a_i|^3,$$

so it follows that $\mu E_1 \geq 1 - \varepsilon$.

Next consider the function

$$f(x) := \sum_{i=1}^n |a_i|^2 |f_i(x)|^2,$$

which satisfies $\|f\|_1 = 1$ and

$$\|f\|_{3/2} \leq \sum_{i=1}^n |a_i|^2 \|f_i^2\|_{3/2} = \sum_{i=1}^n |a_i|^2 \|f_i\|_3^2 \leq M^2.$$

Consider also the set $E_2 := \{x : f(x) > \frac{1}{4}\}$. Since $\int_{E_2^c} f \leq \frac{1}{4}$ it follows

$$\frac{3}{4} \leq \int_{E_2} f \leq \|f\|_{3/2} \mu(E_2)^{1/3} \leq M^2 \mu(E_2)^{1/3}$$

and therefore $\mu E_2 \geq (\frac{3}{4M^2})^3$. Now consider the set

$$E_3 := \Big\{ x \in E_2 : f(x) < 2\Big(\frac{4M^2}{3}\Big)^3 \Big\}.$$

For measure $\mu(E_2 \setminus E_3)$ we have the estimate

$$2\Big(\frac{4M^2}{3}\Big)^3 (\mu E_2 - \mu E_3) \leq \int_{E_2 \setminus E_3} f(x) d\mu(x) \leq \|f\|_1 = 1$$

so that
$$\mu E_3 \geq \frac{1}{2}\left(\frac{3}{4M^2}\right)^3 = 2\varepsilon.$$
Finally note that the set $E := E_1 \cap E_3$ has the following properties
  (i) $\mu E \geq \varepsilon(M) = \frac{1}{4}\left(\frac{3}{4M^2}\right)^3$;
  (ii) $\sum_{i=1}^n |a_i|^3 |f_i(x)|^3 < \frac{256}{27} M^9 \sum_{i=1}^n |a_i|^3$ for all $x \in E$;
  (iii) For $x \in E$ the function $f(x) = \sum_{i=1}^n |a_i|^2 |f_i(x)|^2$ satisfies
$$\frac{1}{4} < f(x) < 2\left(\frac{4M^2}{3}\right) =: \gamma(M).$$

**Step 2.** Define a new measure $\nu$ on $X$ by
$$d\nu(x) = \begin{cases} d\mu(x), & x \in E^c \\ \frac{f(x)\mu(E)}{\int_E f(y)d\mu} d\mu(x), & x \in E. \end{cases}$$
One can easily see that $\nu$ is a probability measure on $X$. Define also functions $g_i(x)$, $1 \leq i \leq n$ by
$$g_i(x) := \begin{cases} f_i(x), & x \in E^c \\ f_i(x)\left(\frac{\int_E f(y)d\mu}{f(x)\mu(E)}\right)^{1/2}, & x \in E. \end{cases}$$
The functions $g_i$ have the following properties:
  (i) $\|g_i\|_{L_2(\nu)} = 1$ for $1 \leq i \leq n$;
  (ii) $\|\sum_{i=1}^n c_i g_i\|_{L_2(\nu)} = \|\sum_{i=1}^n c_i f_i\|_{L_2(\mu)} \leq MR^p \left(\sum_{i=1}^n |c_i|^2\right)^{1/2}$ for all sets of coefficients $\{c_i\}_{i=1}^n$;
  (iii) For all $x \in E$ the following identity takes place
$$g(x) := \sum_{i=1}^n |a_i|^2 |g_i(x)|^2 = \frac{1}{\mu E}\int_E f(y)d\mu(y) =: K^2,$$
  and $\frac{1}{4} \leq K^2 < \gamma(M)$;
  (iv)
$$\sum_{i=1}^n |a_i|^3 |g_i(x)|^3 \leq \beta(M) \sum_{i=1}^n |a_i|^3$$
  for all $x \in E$, where $\beta(M) = 10^5 M^{18}$;
  (v) Finally note that for $x \in E$ and $\omega \in \Omega$ a.s.
$$|\sum_{i=1}^n a_i \xi_i(\omega) g_i(x)| \leq 5M^3 |\sum_{i=1}^n a_i \xi_i(\omega) f_i(x)|.$$

**Step 3.** Note that if there exists a set $F \subset E^m$ such that $\nu^m(F) \geq (\nu E)^m/2$ and for some $\omega_0 \in \Omega$
$$\max_{1 \leq j \leq m}\left(|\sum_{i=1}^n a_i \xi_i(\omega_0) g_i(x_j)|\right) \geq \rho \qquad \text{for } (x_1, \ldots, x_m) \in F,$$
then for the set
$$F_0 = \left\{x \in E : |\sum_{i=1}^n a_i \xi_i(\omega_0) g_i(x)| \geq \rho\right\}$$

we have $(\nu E - \nu F_0)^m \leq (\nu E)^m/2$ so that

$$\mu F_0 \geq C(M)\nu F_0 \geq C(M)\Big[1 - \Big(\frac{1}{2}\Big)^{1/m}\Big]\nu E.$$

Taking into account **(v)** and **(iii)** from step 2 (see also (5)), we get

$$\|\sum_{i=1}^n a_i \xi_i(\omega_0) f_i(x)\|_{\mu,m,\infty} \geq \|\chi_{F_0}(x) \cdot \sum_{i=1}^n a_i \xi_i(\omega_0) f_i(x)\|_{\mu,m,\infty}$$

$$\geq \frac{\rho}{5M^3}\|\chi_{F_0}\|_{\mu,m,\infty} = \frac{\rho}{5M^3}(1 - (1-\mu F_0)^m)$$

$$\geq \frac{\rho}{5M^3}\Big(1 - \Big(1 - C(M)\Big[1 - \Big(\frac{1}{2}\Big)^{\frac{1}{m}}\Big]\nu E\Big)^m\Big).$$

Using the inequality $1 - \left(\frac{1}{2}\right)^{\frac{1}{m}} \geq \frac{1}{2m}$ and **(i)** from step 1 we get

$$\|\sum_{i=1}^n a_i \xi_i(\omega_0) f_i(x)\|_{\mu,m,\infty} \geq \frac{\rho}{5M^3}\Big(1 - \Big(1 - \frac{C(M)\nu E}{2m}\Big)^m\Big) \geq C'(M)\rho.$$

Thus, to prove (9) for $\{f_i\}_1^n$ on $(X,\mu)$ it suffices to prove it for $\{g_i\}_1^n$ on $(X,\nu)$.

Define $F \subset E^m$ by

$$F := \Big\{(x_j)_{j=1}^m \in E^m : \frac{1}{m^2}\sum_{\substack{j,k=1\\j\neq k}}^m |\sum_{i=1}^n |a_i|^2 g_i(x_j)g_i(x_k)|^2 \leq 2\Big(\frac{MR^p}{\varepsilon(M)}\Big)^2 \sum_{i=1}^n |a_i|^4\Big\}.$$

To estimate $\nu F$ notice that

$$\frac{1}{\nu(E)^m}\int_E\cdots\int_E \frac{1}{m^2}\sum_{\substack{j,k=1\\j\neq k}}^m \Big|\sum_{i=1}^n |a_i|^2 g_i(x_j)g_i(x_k)\Big|^2 d\nu(x_1)\ldots d\nu(x_m)$$

$$\leq \frac{1}{(\nu(E)m)^2}\sum_{\substack{j,k=1\\j\neq k}}^m \int_E\int_E \Big|\sum_{i=1}^n |a_i|^2 g_i(x_j)g_i(x_k)\Big|^2 d\nu(x_j)d\nu(x_k)$$

$$\leq \Big[\frac{MR^p}{\nu(E)m}\Big]^2 \sum_{\substack{j,k=1\\j\neq k}}^m \int_E \sum_{i=1}^n |a_i|^4 |g_i(x_j)|^2 d\nu(x_j) \leq \Big[\frac{MR^p}{\varepsilon(M)}\Big]^2 \sum_{i=1}^n |a_i|^4.$$

From Chebyshev's inequality we have $\nu^m(F) \geq \nu(E)^m/2$.

**Step 4.** For $x \in E$ and $\rho > 0$ define

$$E_\rho(x) := \Big\{\omega \in \Omega : \sum_{i=1}^n a_i \xi_i(\omega) g_i(x) > \rho\Big\}.$$

As we have seen in step 3 in order to prove the theorem it suffices to show that there exist some constants $\alpha(M,p) \in (0,1)$, $q' = q'(p) > 0$ and $K_0(M)$ such that for every $(x_j)_{j=1}^m \in F$ and $\rho := \alpha K(2\log P)^{1/2}$ the following estimate takes place

$$(*)\qquad \mathsf{P}\Big\{\Omega \setminus \bigcup_{j=1}^m E_\rho(x_j)\Big\} \leq K_0 P^{-q'}.$$

For fixed $(x_j)_{j=1}^m \in F$ set
$$\eta(\omega) := \sum_{j=1}^m \chi_{E_\rho(x_j)}(\omega).$$

Note that if $\mathsf{P}\big(\bigcup_{j=1}^m E_\rho(x_j)\big) < \kappa$ for some $\kappa$, then by Cauchy-Schwartz inequality we get
$$\mathsf{E}|\eta| \leq (\mathsf{E}|\eta|^2)^{1/2}\mathsf{P}\big(\bigcup_{j=1}^m E_\rho(x_j)\big)^{1/2} < \kappa^{1/2}(\mathsf{E}|\eta|^2)^{1/2}.$$

Thus, the inequality

$$(**) \qquad \mathsf{E}|\eta| \geq (1 - K_0 P^{-q'})(\mathsf{E}|\eta|^2)^{1/2}$$

implies (*) and therefore (9). The aim of the remaining steps is to prove (**).

**Step 5.** In order to prove (**) we shall use a sharper version of the central limit theorem with an estimate for the error term. We use the following result due to Rotar' [13] (or see Corollary 17.2 in [1]).

**Proposition 1.** *Let $\{X_i\}_{i=1}^h$ be a set of random vectors in $\mathbb{R}^d$ such that $\mathsf{E}X_i = 0$, $1 \leq i \leq h$ then*
$$\sup_{A \in \mathcal{C}} |Q_h(A) - \Phi_{0,V}(A)| \leq K_1(d) h^{-1/2} \rho_3 \lambda^{-3/2},$$

*where $K_1(d) < \infty$ is a constant, $\mathcal{C}$ denotes the class of all Borel convex subsets of $\mathbb{R}^d$,*
$$\rho_3 := h^{-1} \sum_{i=1}^h \mathsf{E}|X_i|^3,$$

*$\lambda$ is the smallest eigenvalue of the matrix $V = h^{-1} \sum_1^h \mathrm{cov}(X_i)$ ( recall that the covariance matrix of a random vector $Y = (y_1, \ldots, y_d)$ such that $\mathsf{E}Y = 0$ defined by $\mathrm{cov}(Y) := \{\mathsf{E}(y_j y_k)\}_{j,k=1}^d$), $Q(A)$ is the probability that $h^{-1/2} \sum_{i=1}^h X_i$ belongs to a convex set $A$ and, finally, $\Phi_{0,V}$ denotes the normal distribution with the density*
$$\phi_{0,V}(Y) := (2\pi)^{-d/2} (\det V)^{-1/2} \exp\big\{-\frac{1}{2}(Y, V^{-1}Y)\big\}, \qquad Y \in \mathbb{R}^d.$$

We shall apply Proposition 1 twice: for one- and two-dimensional cases.

For fixed $x \in E$ let
$$X_i(\omega) := a_i \xi_i(\omega) g_i(x) \qquad 1 \leq i \leq n.$$

Then
$$\rho_3 = \frac{1}{n} \sum_{i=1}^n \mathsf{E}|X_i|^3 = \frac{1}{n} \sum_{i=1}^n |a_i|^3 \mathsf{E}|\xi_i|^3 |g_i(x)|^3 \leq$$
$$\leq \frac{M^3}{n} \sum_{i=1}^n |a_i|^3 |g_i(x)|^3 \leq \frac{M^3 \beta(M)}{n} \sum_{i=1}^n |a_i|^3,$$
$$\lambda = V = \frac{1}{n} \sum_{i=1}^n \mathrm{cov}(X_i) = \frac{1}{n} \sum_{i=1}^n |a_i|^2 \mathsf{E}|\xi_i|^2 |g_i(x)|^2 = \frac{K^2}{n}.$$

So Preposition 1 implies

$$\left|\mathsf{P}(E_\rho(x)) - \frac{1}{K}\left(\frac{n}{2\pi}\right)^{1/2}\int_{\frac{\rho}{\sqrt{n}}}^\infty e^{-\frac{y^2 n}{2K^2}}dy\right| \leq K_1\frac{M^3\beta(M)}{K^3}\sum_{i=1}^n |a_i|^3.$$

By a change of variable in the integral we get

$$\left|\mathsf{P}(E_\rho(x)) - \frac{1}{(2\pi)^{1/2}K}\int_\rho^\infty e^{-\frac{y^2}{2K^2}}dy\right| \leq K_2(M)\sum_{i=1}^n |a_i|^3 \leq K_2\Big(\sum_{i=1}^n |a_i|^4\Big)^{\frac{1}{2}} = \frac{K_2}{R^{1/2}}.$$

It is well-known that

$$\int_z^\infty e^{-t^2/2}dt \asymp \frac{1}{z}e^{-z^2/2}, \qquad z > 1.$$

Therefore, when $\alpha(M,p)$ satisfies $0 < \alpha^2 < 1/2$ and $R > R_0(M)$ we can neglect the error term in the application of the central limit theorem, so we have

$$\mathsf{P}(E_\rho(x)) \asymp K\rho^{-1}e^{-\frac{\rho^2}{2K^2}}$$

which implies

$$\mathsf{E}|\eta| := \sum_{j=1}^m \mathsf{P}(E_\rho(x_j)) \asymp mK\rho^{-1}e^{-\frac{\rho^2}{2K^2}}.$$

Note here, that by taking if necessary $K_0(M)$ large enough we can neglect the case $R < R_0(M)$.

Note that

$$\mathsf{E}|\eta|^2 = \sum_{j=1}^m \mathsf{P}(E_\rho(x_j)) + \sum_{\substack{j,k=1 \\ j\neq k}}^m \mathsf{P}(E_\rho(x_j) \cap E_\rho(x_k)) =$$

$$= \mathsf{E}|\eta| + \sum_{\substack{j,k=1 \\ j\neq k}}^m \mathsf{P}(E_\rho(x_j) \cap E_\rho(x_k))$$

and $\mathsf{E}|\eta| \leq K_3 m^{-1/2}(\mathsf{E}|\eta|)^2 \leq K_3 P^{-1/2}(\mathsf{E}|\eta|)^2$, where $K_3 = K_3(M) > 0$. So to prove (**) it is enough to show that

$$(***) \qquad \sum_{\substack{j,k=1 \\ j\neq k}}^m \mathsf{P}\big(E_\rho(x_j) \cap E_\rho(x_k)\big) \leq (1 + K_4 P^{-q'})(\mathsf{E}|\eta|)^2$$

for some constants $K_4 = K_4(M,p) > 0$, $q' = q'(M) > 0$, $\alpha(M,p) > 0$.

**Step 6.** Let us split the index set $\{(i,j) : 1 \leq i \neq j \leq n\}$ into two sets. Let

$$\sigma_1 = \Big\{(j,k) : 1 \leq j \neq k \leq m, |\sum_{i=1}^n |a_i|^2 g_i(x_j)g_i(x_k)| < \frac{1}{8}\Big\}.$$

Since $(x_j)_{j=1}^n \in F$ (see Step 3) it follows that

$$|\sigma_1^c| \leq 8^2 \sum_{\substack{j,k=1 \\ j\neq k}}^m |\sum_{i=1}^n |a_i|^2 g_i(x_j)g_i(x_k)|^2 \leq 128\Big(\frac{mMR^p}{\varepsilon(M)}\Big)^2 \sum_{i=1}^n |a_i|^4 = \frac{128}{R}\Big(\frac{mMR^p}{\varepsilon(M)}\Big)^2.$$

Thus, whenever $\alpha^2(M,p) < 1/2 - p$ we have

$$\sum_{(j,k)\in\sigma_1^c} \mathsf{P}\big(E_\rho(x_j) \cap E_\rho(x_k)\big) \leq \frac{128}{R}\Big(\frac{mMR^p}{\varepsilon(M)}\Big)^2 \rho^{-1} e^{-\frac{\rho^2}{2K^2}} < \frac{K_4(M)}{R^{1/2-p}}(\mathsf{E}|\eta|)^2.$$

**Step 7.** For fixed pair $s = (j,k) \in \sigma_1$ consider a set of 2-dimensional random vectors defined by

$$X_i^s(\omega) := (a_i \xi_i(\omega) g_i(x_j), a_i \xi_i(\omega) g_i(x_k)); \qquad 1 \leq i \leq n.$$

To estimate the error term in the central limit theorem for these random vectors, notice that

$$\rho_3^s := \frac{1}{n}\sum_{i=1}^n |a_i|^3 \mathsf{E}|\xi_i|^3 \big(|g_i(x_j)|^2 + |g_i(x_k)|^2\big)^{3/2} \leq \frac{8M^3 \beta(M)}{n} \sum_{i=1}^n |a_i|^3,$$

$$V^s = \frac{1}{n}\begin{pmatrix} \sum_{i=1}^n |a_i|^2 |g_i(x_j)|^2 & \sum_{i=1}^n |a_i|^2 g_i(x_j) g_i(x_k) \\ \sum_{i=1}^n |a_i|^2 g_i(x_j) g_i(x_k) & \sum_{i=1}^n |a_i|^2 |g_i(x_k)|^2 \end{pmatrix}$$

$$= \frac{1}{n}\begin{pmatrix} K^2 & \sum_{i=1}^n |a_i|^2 g_i(x_j) g_i(x_k) \\ \sum_{i=1}^n |a_i|^2 g_i(x_j) g_i(x_k) & K^2 \end{pmatrix}.$$

Hence,

$$\det V^s = \frac{1}{n^2}\Big(K^4 - \big|\sum_{i=1}^n |a_i|^2 g_i(x_j) g_i(x_k)\big|^2\Big) > \frac{1}{n^2}\Big(\frac{1}{16} - \frac{1}{64}\Big) = \frac{3}{64n^2},$$

and

$$\operatorname{trace} V^s = \frac{2K^2}{n}.$$

Note that the matrix $V^s$ is positive so both its eigenvalues are positive. Let $\lambda_2 \geq \lambda_1 > 0$ be the eigenvalues, taking into account that $\lambda_1 + \lambda_2 = \operatorname{trace} V^s = 2K^2/n$, we get

$$\frac{3}{64n^2} < \det V^s = \lambda_2 \lambda_1 < \frac{2K^2}{n}\lambda_1,$$

hence,

$$\lambda_1 > \frac{3}{128nK^2}.$$

So the central limit theorem (Proposition 1) for $X_i^s$ gives

$$\Big|\mathsf{P}(E_\rho(x_j) \cap E_\rho(x_k)) - \frac{1}{2\pi(\det V^s)^{1/2}} \int_{\frac{\rho}{\sqrt{n}}}^\infty \int_{\frac{\rho}{\sqrt{n}}}^\infty e^{-\frac{1}{2}(Y,(V^s)^{-1}Y)} dy_1 dy_2 \Big| \leq$$

$$\leq K_1(2) n^{-1/2} \rho_3^s \lambda_1^{-3/2} < K_5 \sum_{i=1}^n |a_i|^3,$$

for a constant $K_5 = K_5(M) < \infty$. Taking into account $\sum_{i=1}^n |a_i|^3 \leq R^{-1/2}$, we get

$$\sum_{s=(j,k)\in\sigma_1} \mathsf{P}(E_\rho(x_j) \cap E_\rho(x_k)) <$$

$$< K_5 \frac{m^2}{R^{1/2}} + \int_\rho^\infty \int_\rho^\infty \sum_{s\in\sigma_1} \frac{1}{2\pi(\det nV^s)^{1/2}} e^{-\frac{1}{2}(Y,(nV^s)^{-1}Y)} dy_1 dy_2.$$

If we choose $\alpha(M) < 1/5$ the error term $K_5 m^2 R^{-1/2} \leq K_6(M) R^{-1/4} (\mathsf{E}|\eta|)^2$, thus, to prove (***) it remains to estimate the integral term

$$\int_\rho^\infty \int_\rho^\infty \sum_{s\in\sigma_1} \frac{1}{2\pi(\det nV^s)^{1/2}} e^{-\frac{1}{2}(Y,(nV^s)^{-1}Y)} dy_1 dy_2.$$

We shall compare it with the expression

$$\sum_{(j,k)\in\sigma_1} \mathsf{P}(E_\rho(x_j)) \cdot \mathsf{P}(E_\rho(x_k)) = \int_\rho^\infty \int_\rho^\infty \sum_{s\in\sigma_1} \frac{1}{2\pi K^2} e^{-\frac{1}{2K^2}(y_1^2+y_2^2)} dy_1 dy_2 + w$$

$$= \frac{|\sigma_1|}{2\pi K^2} \int_\rho^\infty \int_\rho^\infty e^{-\frac{1}{2K^2}(y_1^2+y_2^2)} dy_1 dy_2 + w,$$

where $w \leq 2K_2 m^2 R^{-1/2}$ (see Step 5), as before we can ensure that $w = o(R^{-1/4})(\mathsf{E}|\eta|)^2$. In order to compare the two integral expressions, notice that

$$(nV^s)^{-1} = \frac{1}{\det(nV^s)} \begin{pmatrix} K^2 & -\sum_{i=1}^n |a_i|^2 g_i(x_j) g_i(x_k) \\ -\sum_{i=1}^n |a_i|^2 g_i(x_j) g_i(x_k) & K^2 \end{pmatrix}.$$

Let

$$c_s := \sum_{i=1}^n |a_i|^2 g_i(x_j) g_i(x_k) \qquad \text{for} \quad s = (i,j) \in \sigma_1,$$

then

$$(nV^s)^{-1} = \begin{pmatrix} \frac{K^2}{K^4-|c_s|^2} & -\frac{c_s}{K^4-|c_s|^2} \\ -\frac{c_s}{K^4-|c_s|^2} & \frac{K^2}{K^4-|c_s|^2} \end{pmatrix}.$$

Now we have

$$\int_\rho^\infty \int_\rho^\infty \sum_{s\in\sigma_1} \frac{1}{2\pi(\det nV^s)^{1/2}} e^{-\frac{1}{2}(Y,(nV^s)^{-1}Y)} dy_1 dy_2 =$$

$$= \int_\rho^\infty \int_\rho^\infty \sum_{s\in\sigma_1} \frac{1}{2\pi(K^4-|c_s|^2)^{1/2}} \exp\left(-\frac{K^2(y_1^2+y_2^2)}{2(K^4-|c_s|^2)} + \frac{c_s y_1 y_2}{K^4-|c_s|^2}\right) dy_1 dy_2.$$

Let also

$$a_s := \frac{K^2}{K^4-|c_s|^2}; \qquad b_s := \frac{c_s}{K^4-|c_s|^2}, \qquad s \in \sigma_1.$$

Notice, that for any $L > 1$ the following inequality holds

$$\int_{L\rho}^{\infty} \int_{\rho}^{\infty} \sum_{s \in \sigma_1} \frac{1}{2\pi(K^4 - |c_s|^2)^{1/2}} e^{\{-\frac{a_s}{2}(y_1^2 + y_2^2) + b_s y_1 y_2\}} dy_1 dy_2 \leq$$

$$\leq 4 \int_{L\rho}^{\infty} \int_{\rho}^{\infty} \sum_{s \in \sigma_1} e^{-\frac{a_s - b_s}{2}(y_1^2 + y_2^2)} dy_1 dy_2 \asymp \frac{1}{L\rho^2} \sum_{s \in \sigma_1} e^{-\frac{a_s - b_s}{2}(L^2 + 1)\rho^2}.$$

Since for all $s \in \sigma_1$ we have

$$a_s - b_s = \frac{K^2 - c_s}{K^4 - |c_s|^2} = \frac{1}{K^2 + c_s} \geq \frac{1}{\gamma(M) + \frac{1}{8}} =: \gamma'(M) > 0,$$

it follows that

$$\int_{L\rho}^{\infty} \int_{\rho}^{\infty} \sum_{s \in \sigma_1} \frac{1}{2\pi(\det(nV^s))^{1/2}} e^{\{-\frac{a_s}{2}(y_1^2 + y_2^2) + b_s y_1 y_2\}} dy_1 dy_2 \leq \frac{K_7 m^2}{L\rho^2} e^{-\frac{1}{2}\gamma'(L^2 + 1)\rho^2}.$$

Set $L^2 + 1 := \frac{4}{\gamma'(M) K^2}$ and get

$$\frac{K_7 m^2}{L\rho^2} e^{-\frac{1}{2}\gamma'(M)(L^2+1)\rho^2} \leq K_8(M) \frac{(\mathsf{E}|\eta|)^2}{P^{\alpha^2(M,p)}}$$

with a constant $K_8(M) < \infty$. Therefore, there exists a constant $K_9(M) < \infty$ such that whenever $\alpha^2(M, p) < 1/5$ the following inequality holds

$$\sum_{(j,k) \in \sigma_1} \mathsf{P}(E_\rho(x_j)) \cap \mathsf{P}(E_\rho(x_k)) \leq K_9(M) \frac{(\mathsf{E}|\eta|)^2}{P^{\alpha^2}} +$$

$$+ \int_{\rho}^{L\rho} \int_{\rho}^{L\rho} \sum_{s \in \sigma_1} \frac{\exp\left\{-\frac{a_s}{2}(y_1^2 + y_2^2) + b_s y_1 y_2\right\}}{2\pi(\det(nV^s))^{1/2}} dy_1 dy_2.$$

**Step 8.** To finish the proof of the Theorem it remains to compare the expression

$$A := \sum_{s \in \sigma_1} \frac{1}{2\pi(\det(nV^s))^{1/2}} e^{\{-\frac{a_s}{2}(y_1^2 + y_2^2) + b_s y_1 y_2\}}$$

with the expression

$$B := \frac{1}{2\pi K^2} |\sigma_1| e^{-\frac{1}{2K^2}(y_1^2 + y_2^2)}$$

in the range $\rho \leq y_1, y_2 \leq L\rho$. We are going to show that $A \leq B(1 + K_{10} R^{-q'})$ pointwise in that range with some constants $K_{10}(M) < \infty$, $q'(p) > 0$. In fact, assume for a moment we have shown it, then integrate this inequality over the domain $\rho \leq y_1, y_2 \leq L\rho$ and get

$$\int_{\rho}^{L\rho} \int_{\rho}^{L\rho} \sum_{s \in \sigma_1} \frac{1}{2\pi(\det(nV^s))^{1/2}} e^{\{-\frac{a_s}{2}(y_1^2 + y_2^2) + b_s y_1 y_2\}} dy_1 dy_2 <$$

$$< (1 + \frac{K_{10}}{R^{q'}}) \int_{\rho}^{\infty} \int_{\rho}^{\infty} \frac{|\sigma_1|}{2\pi K^2} e^{-\frac{y_1^2 + y_2^2}{2K^2}} dy_1 dy_2 \leq (1 + \frac{K_{10}}{R^{q'}})((\mathsf{E}|\eta|)^2 + w),$$

where the error term $w = o(R^{-1/4})(\mathsf{E}|\eta|)^2$ (see step 7). This finally implies (***) and, thus, the theorem statement.

To prove this inequality split the index set $\sigma_1$ into subsets
$$\sigma_r := \{s \in \sigma_1 : 2^{-r} \leq |c_s| < 2^{-r+1}\}, \qquad r = 4, 5, \ldots.$$
Clearly, $\sigma_1 = \bigcup_{r \geq 4} \sigma_r$. We can estimate $|\sigma_r|$ as follows (see Step 3)
$$\frac{|\sigma_r|}{2^{2r}} \leq \sum_{s \in \sigma_1} |c_s|^2 \leq 2m^2 \Big(\frac{MR^p}{\varepsilon(M)}\Big)^2 \sum_{i=1}^{n} |a_i|^4$$
so that
$$|\sigma_r| \leq \min\Big(m^2, 2^{2r+1}\Big(\frac{MR^p}{\varepsilon(M)}\Big)^2 \frac{m^2}{R}\Big).$$
For $s \in \sigma_r$ and $\rho \leq y_1, y_2 \leq L\rho$ we have
$$\Big|(a_s - \frac{1}{K^2})(y_1^2 + y_2^2) - 2b_s y_1 y_2\Big| = \Big|\frac{(c_s/K)^2(y_1^2 + y_2^2) - 2c_s y_1 y_2}{K^4 - c_s^2}\Big| \leq$$
$$\leq 25\big(2^{-2r+2}(y_1^2 + y_2^2) + 2^{-r+2} y_1 y_2\big) \leq 200 \cdot 2^{-r} L^2 \rho^2.$$
Moreover, if $s \in \sigma_r$ with $r \geq 4$, then
$$\det(nV^s) = K^4 - c_s^2 \geq K^4 - 2^{-2r+2} \geq \frac{K^4}{(1 + K_{11} 2^{-r})^2},$$
where $K_{11}$ is an absolute constant. Now we can say that
$$A \leq \frac{B \cdot S}{|\sigma_1|},$$
where
$$S := \sum_{r=4}^{\infty} |\sigma_r|(1 + K_{11} 2^{-r}) e^{\frac{200}{2^r} L^2 \rho^2} = \sum_{r=4}^{[\frac{1-2p}{4} \log_2 R]} + \sum_{[\frac{1-2p}{4} \log_2 R]+1}^{\infty} = S_1 + S_2.$$

We can estimate $S_1$ as follows
$$S_1 \leq \frac{1-2p}{4} \log_2 R \cdot 2R^{\frac{1-2p}{2}} \Big(\frac{MR^p}{\varepsilon(M)}\Big)^2 \frac{m^2}{R}\Big(1 + \frac{K_{11}}{16}\Big) e^{\frac{200}{16} L^2 \rho^2}.$$
Notice, that $e^{\frac{200}{16} L^2 \rho^2} = P^{\frac{200}{8} \alpha^2 K^2 L^2}$. Choose $\alpha(M, p)$ such that
$$\frac{200}{8} \alpha^2 K^2 L^2 < \frac{1-2p}{20}.$$
This condition on $\alpha(M, p)$ is compatible with the previously imposed ones ($\alpha^2 < \min\{\frac{1}{2} - p, \frac{1}{25}\}$). So we get
$$S_1 \leq K_{12} |\sigma_1| R^{-\frac{1-2p}{3}}$$
with a constant $K_{12} = K_{12}(M, p) < \infty$. Further,
$$S_2 \leq \sum_{r > [\frac{1-2p}{4} \log_2 R]} |\sigma_r|(1 + K_{11} R^{\frac{-(1-2p)}{4}})(1 + K_{13} R^{\frac{-(1-2p)}{5}}) \leq |\sigma_1|(1 + K_{14} R^{\frac{-(1-2p)}{5}})$$

with some constants $K_{13} = K_{13}(M,p)$ and $K_{14} = K_{14}(M,p) < \infty$. Thus, we have
$$|S| \leq \big(1 + (K_{12} + K_{14})R^{-(1-2p)/5}\big)|\sigma_1|$$
which implies $A \leq B\big(1 + K_{10}R^{-(1-2p)/5}\big)$ and completes the proof of Theorem 1.

## 3. THE UPPER ESTIMATE.

Let us show now that with some restrictions on $\{\xi_k\}_1^n$ and $m \leq n$ the estimate (10) is precise in sense of order. The following theorem states this explicitly.

**Theorem 2.** *Let $\xi_k$ be independent variables for which the following exponential estimate takes place*

$$\tag{11} \mathsf{P}\Big\{\big|\sum_{k=1}^n c_k \xi_k\big| > t\big(\sum_{k=1}^n c_k^2\big)^{1/2}\Big\} \leq C_4 e^{-t^2 C_5}$$

*for all sets of coefficients $\{c_k\}_1^n$ with some absolute positive constants $C_4$, $C_5$. Then*

$$\tag{12} \mathsf{E}\|\sum_{k=1}^n \xi_k f_k\|_{m,\infty} \leq C_6 \big\|\big(\sum_{k=1}^n |f_k|^2\big)^{1/2}\big\|_{m,\infty} \cdot \sqrt{1 + \log m}$$

*for all function systems $\{f_k\}_1^n \subset L_1(X)$ and all $m \geq 1$ with an absolute constant $C_6 > 0$. And since for all bounded functions $f \in L_\infty$ the integral-uniform norm $\|f\|_{m,\infty} \leq \|f\|_\infty$ for bounded functions (12) implies*

$$\mathsf{E}\|\sum_{k=1}^n \xi_k f_k\|_{m,\infty} \leq C_6 \Big(\sum_{k=1}^n \|f_k\|_\infty^2\Big)^{1/2} \cdot \sqrt{1 + \log m}.$$

**Proof.** For $x \in X$ define
$$\eta_x(\omega) := \sum_{k=1}^n f_k(x)\xi_k(\omega);$$
$$\mu_t(x_1,\ldots,x_m) := \mathsf{P}\Big\{\max_{1 \leq h \leq m}(|\eta_{x_h}|) > t \max_{1 \leq h \leq m}\big(\sum_{k=1}^n |f_k(x_h)|^2\big)^{1/2}\Big\}.$$

Using the exponential estimate (11), we get that
$$\mu_t(x_1,\ldots,x_m) \leq \sum_{h=1}^m \mathsf{P}\Big\{|\eta_{x_h}| > t \max_{1 \leq h \leq m}\big(\sum_{k=1}^n |f_k(x_h)|^2\big)^{1/2}\Big\} \leq mC_4 e^{-C_5 t^2}.$$

Note that for $t_0 = \big(\frac{3}{C_5}\log m\big)^{1/2}$ we have $\mu_{t_0}(x_1,\ldots,x_m) \leq C_4 m^{-2}$. Now we can estimate

$$\mathsf{E}\max_{1\leq h \leq m}|\eta_{x_h}| \leq \Big(\max_{1\leq h\leq m}\sum_{k=1}^n |f_k(x_h)|^2\Big)^{1/2}\Big(\sqrt{\frac{3}{C_5}\log m} + \sum_{t=t_0}^\infty (t+1)\mu_t(x_1,\ldots,x_m)\Big)$$

$$\leq \Big(\max_{1\leq h\leq m}\sum_{k=1}^n |f_k(x_h)|^2\Big)^{1/2}\Big(\sqrt{\frac{3}{C_5}\log m} + C_4 m \sum_{t=t_0}^\infty (t+1)e^{-C_5 t^2}\Big)$$

$$\leq C \max_{1\leq h\leq m}\Big(\sum_{k=1}^n |f_k(x_h)|^2\Big)^{1/2}\sqrt{1 + \log m}.$$

Integrating the last inequality with respect to $x_1, \ldots, x_m$, we get (see (2))

$$\mathsf{E}\|\sum_{k=1}^{n}\xi_k f_k\|_{m,\infty} \leq \int_X \ldots \int_X \mathsf{E}\max_{1\leq h\leq m}|\eta_{x_h}|d\mu(x_1)\ldots d\mu(x_m)$$

$$\leq C\left\|\left(\sum_{k=1}^{n}|f_k|^2(1+\log m)\right)^{1/2}\right\|_{m,\infty}.$$

This completes the proof.

**Corollary 2.** *For uniformly bounded functions $\{f_k\}_1^n \subset L_\infty(X)$ with $\|f_k\|_\infty \leq M$ and independent random variables $\{\xi_k\}_1^n$, satisfying the exponential estimate (11), Theorem 2 implies*

$$\mathsf{E}\|\sum_{k=1}^{n}a_k\xi_k f_k\|_{m,\infty} \leq MC_6\left(\sum_{k=1}^{n}|a_k|^2\right)^{1/2}\cdot\sqrt{1+\log m}.$$

Hence, whenever $m \leq n$ and $m = O\bigl(R(\{a_k\}_1^n)\bigr)$ (see (7)), then the inequality (10) from Theorem 1 is precise in sense of order for all uniformly bounded function systems $\{f_k\}_1^n$ and independent random variables $\{\xi_k\}_1^n$, satisfying the exponential estimate (11). In particular, it is true for trigonometric polynomials with random coefficients.

**Remark.** If we take a sequence of (multivariate) trigonometric polynomials of order at most $n$ as the functions $\{f_k\}_1^n$ and apply Theorem 2 with the parameter $m = n$ then, taking into account (14) (see below), we get the well-known upper estimate for the expectation of the *uniform* norm of a random trigonometric polynomial, e.g. exposed in J.-P. Kahane's book (see Th. 3 Ch. 6 [4]).

## 4. SOME PROPERTIES OF THE INTEGRAL-UNIFORM NORMS AND APPLICATION.

The following Theorem compares the integral-uniform norm of an integrable function $f \in L_1(X)$ with its average over an arbitrary subset of $X$.

**Theorem 3.** *For each $f \in L_1(X)$ ($(X,\mu)$ is a probability space) and arbitrary measurable $\Delta \subset X$ ($\mu\Delta \equiv |\Delta| > 0$) the following inequality holds*

(13) $$\|f\|_{m,\infty} \geq \bigl(1-(1-|\Delta|)^m\bigr)\cdot\frac{1}{|\Delta|}\int_\Delta|f|.$$

**Proof.** Clearly, it suffices to prove (13) for the case when $\text{supp} f \subset \Delta$. Using the formula (3) we get

$$\|f\|_{k+1,\infty} - \|f\|_{k,\infty} = \int_0^\infty \bigl(1-(1-\lambda_f(t))^{k+1}\bigr)dt - \int_0^\infty \bigl(1-(1-\lambda_f(t))^k\bigr)dt$$

$$= \int_0^\infty \lambda_f(t)\bigl(1-\lambda_f(t)\bigr)^k dt$$

$$\geq (1-|\Delta|)^k\int_0^\infty \lambda_f(t)dt = (1-|\Delta|)^k\|f\|_1.$$

Sum this inequality up from $k = 1$ to $k = m-1$ and get

$$\|f\|_{m,\infty} - \|f\|_1 \geq \|f\|_1\sum_{k=1}^{m-1}(1-|\Delta|)^k,$$

which implies
$$\|f\|_{m,\infty} \geq \|f\|_1 \frac{1 - (1 - |\Delta|)^m}{|\Delta|}.$$
To complete the proof notice that $\|f\|_1 = \int_\Delta |f|$ since $\operatorname{supp} f \subset \Delta$.

For trigonometric polynomials of order at most $n$ the identity (5) implies
$$\|P_n\|_{n,\infty} \asymp \|P_n\|_\infty. \tag{14}$$
In fact, for a set $E := \{x \in [0, 2\pi] : |P_n(x)| \geq \|P_n\|_\infty/2\}$ the Bernstein inequality implies $\mu E \geq 1/n$, evaluating $\|\chi_E\|_{n,\infty}$ from (5), we get (14). If $n \geq m$ then (4) and (14) for trigonometric polynomials of order at most $n$ imply
$$\|P_n\|_\infty \leq C \frac{n}{m} \|P_n\|_{m,\infty},$$
where $C > 0$ is an absolute constant. For the Fejér kernels this inequality is precise in sense of order, in fact, when $n \geq m$ one can prove that
$$\|K_n\|_{m,\infty} \asymp m;$$
$$\|D_n\|_{m,\infty} \asymp m(1 + \log \frac{n}{m}),$$
where $K_n$ is the Fejér kernel and $D_n$ is the Dirichlet kernel.

For the integral-uniform norm as for any shift invariant norm (e.g. see [2]) the following analog of the Bernstein inequality takes place.

**Proposition 2.**[1] *For the integral-uniform norm of the derivative of trigonometric polynomial $P_n$ of order at most $n$ the following inequality holds*
$$\|P_n^{(r)}\|_{m,\infty} \leq n^r \|P_n\|_{m,\infty}, \qquad r = 1, 2, \ldots. \tag{15}$$

**The idea of the proof.** Use the M. Riesz Interpolation Formula [12] (or see Ch. 2.4 [11]) for derivative of a trigonometric polynomial of order at most $n$:
$$P_n'(x) = \sum_{k=1}^{2n} (-1)^{k+1} \lambda_k P_n(x + x_k),$$
$$\text{where} \qquad \lambda_k := \frac{1}{4n \sin^2(x_n/2)}; \qquad x_n := \frac{2k-1}{2n}\pi.$$
And notice that $\sum_{k=1}^{2n} \lambda_k = n.\square$

It is well-known that $L_\infty$-norm of trigonometric polynomials of order at most $n$ is equivalent to its $L_\infty$-norm taken over the uniform net $\{\frac{s}{4n} 2\pi\}_{s=1}^{4n}$, precisely
$$\|P_n\|_\infty \asymp \max_{1 \leq s \leq 4n} \left(|P_n(\frac{s}{4n} 2\pi)|\right). \tag{16}$$

This fact easily follows from the classical Bernstein inequality for the uniform norm. Using Proposition 2 one can prove an analog of (16) for the integral-uniform norm. For

---
[1] See [2] for more general cases of Bernstein-type inequalities.

a vector $\mathbf{x} = (x_k)_{k=1}^N \in \mathbb{R}^N$ define $\|\mathbf{x}\|_{m,\infty}$ by

$$\|\mathbf{x}\|_{m,\infty} := \frac{1}{N^m} \sum_{k_1=1}^N \cdots \sum_{k_m=1}^N \max_{1 \leq j \leq m} |x_{k_j}|.$$

For trigonometric polynomials we have the following

**Theorem 4.** *There exist positive constants $C_9$, $C_{10}$ such that for all trigonometric polynomials $P_n$ of order at most $n$ the following inequalities hold*

(17) $$C_9 \|\mathbf{p_n}\|_{m,\infty} \leq \|P_n\|_{m,\infty} \leq C_{10} \|\mathbf{p_n}\|_{m,\infty},$$

where $\mathbf{p_n} := \bigl(P_n(t_k)\bigr)_{k=1}^{8n}$, $t_k := \frac{k-1}{8n} 2\pi$.

**Proof.** Let $\Delta_k := [t_k, t_{k+1})$, and $\delta := 2\pi/(8n)$. The family of semi-intervals $\{\Delta_k\}_{k=1}^{8n}$ splits the circle $[0, 2\pi]$, so for each $x \in [0, 2\pi)$ there exists a unique $k(x)$ such that $x \in \Delta_{k(x)}$. Thus, for any net $x_1, \ldots, x_m$ we have

$$\Bigl| \max_{1 \leq j \leq m} \bigl(|P_n(x_j)|\bigr) - \max_{1 \leq j \leq m} \bigl(|P_n(t_{k(x_j)})|\bigr) \Bigr| \leq \max_{1 \leq j \leq m} \bigl(|P_n(x_j) - P_n(t_{k(x_j)})|\bigr)$$

$$\leq \max_{1 \leq j \leq m} \Bigl( \int_{\Delta_{k(x_j)}} |P_n'| \Bigr).$$

Integrating this inequality over $x_1, \ldots, x_m$ we get

$$(2\pi)^m \bigl| \|P_n\|_{m,\infty} - \|\mathbf{p_n}\|_{m,\infty} \bigr| \leq \delta^m \sum_{j_1,\ldots,j_m} \max_{1 \leq s \leq m} \int_{\Delta_{j_s}} |P_n'| \leq$$

$$\leq \delta \sum_{j_1,\ldots,j_m} \int_{\Delta_1} \cdots \int_{\Delta_m} \max_{1 \leq s \leq m} |P_n'(x_s)| dx_1 \ldots dx_m = \delta (2\pi)^m \|P_n'\|_{m,\infty}.$$

Applying (15) to estimate the righthand-side we get

$$\bigl| \|P_n\|_{m,\infty} - \|\mathbf{p_n}\|_{m,\infty} \bigr| \leq \delta n \|P_n\|_{m,\infty} < \frac{2\pi}{8} \|P_n\|_{m,\infty}.$$

Since $\pi/4 < 1$, it implies (17) and completes the proof.

Now we give an application of Theorem 1 demonstrating the potential utility of the family of integral-uniform norms. In [10] S. Montgomery-Smith and E.M. Semenov reduced a certain problem from functional analysis to the following

**Problem.** *Let $\{f_i\}_{i=1}^n$ be a system of functions defined on a measure space $(X, \mu)$, $\mu X = 1$ and $\|f_i\|_1 = 1$. The question is if there exists a sequence of signs $\{\theta_i\}_{i=1}^n$, $\theta_i = \pm 1$ such that for every $k = 1, \ldots, n$ the following estimate takes place*

(18) $$\sup_{\substack{e \subset X \\ \mu e = 2^{-k}}} 2^k \int_e |\sum_{i=1}^n \theta_i f_i(x)| d\mu(x) \geq c_0 \sqrt{nk},$$

*where $c_0$ is an absolute positive constant.*

Using Theorem 1 we prove the following theorem, which partially solves the Problem.

**Theorem 5.** *For a set of functions $\{f_i\}_{i=1}^n$ on $(X, \mu)$, $\mu X = 1$ such that $\|f_i\|_1 = 1$ and $\|f_i\|_3 \leq M$, there exists a sequence of signs $\{\theta_i\}_{i=1}^n$, $\theta_i = \pm 1$, such that for every $k = 1, \ldots, \log n$ the inequality (18) takes place with an absolute constant $c_0 > 0$.*

**Proof.** Fix some $\delta \in (0, 1/2)$. Assume first that for the function system $\{f_i\}_{i=1}^n$ there exists a set of signs $\{\theta_i\}_{i=1}^n$, $\theta_i = \pm 1$ such that

$$\|\sum_{i=1}^n \theta_i f_i\|_1 \geq n^{\frac{1}{2}+\delta}.$$

Then it obviously implies the assertion of the Theorem. In fact, let $k < n^\delta$ then

$$\sup_{\substack{e \subset X \\ \mu e = 2^{-k}}} 2^k \int_e |\sum_{i=1}^n \theta_i f_i(x)| d\mu(x) \geq \|\sum_{i=1}^n \theta_i f_i\|_1 \geq n^{1/2+\delta}.$$

Now assume the opposite, i.e. that

(19) $$\|\sum_{i=1}^n \theta_i f_i\|_1 \leq n^{\frac{1}{2}+\delta}.$$

for all sets of signs $\{\theta_i\}_{i=1}^n$, $\theta_i = \pm 1$. We are going to show that (19) and the boundness of functions $f_i$ in $L_3$ imply **(b')** from Section 2 with some $p < 1/2$. We need the following geometrical

**Lemma.** *Let $\{w_i\}_{i=1}^n$ be a set of vectors in a linear space with a norm $\|\cdot\|$ (or a seminorm) such that $\|w_i\| = 1$ and*

(20) $$\|\sum_{i=1}^n \theta_i w_i\| \leq C_{11} n^{\frac{1}{2}+\beta}$$

*for all sets of signs $\{\theta_i\}_{i=1}^n$, $\theta_i = \pm 1$ with some constants $\beta \in [0, 1/2)$, $C_{11} > 0$. Then*

(20') $$\|\sum_{i=1}^n a_i w_i\| \leq C_{12}(\beta, C_{11}) n^{\frac{1}{4}+\frac{\beta}{2}} \Big(\sum_{i=1}^n a_i^2\Big)^{\frac{1}{2}}$$

*for all sets of coefficients $\{a_i\}_1^n$. One cannot improve this estimate in the sense that there exist a norm $\|\cdot\|$, vectors $\{w_i\}_1^n$ and coefficients $\{a_i\}_1^n$ such that (20) holds and (20') is precise in sense of order.*[2]

Let us first finish the proof of the theorem. Notice that (19), Hölder's inequality and the triangle inequality for $\|\cdot\|_3$-norm imply for the function $F^\theta := \sum_1^n \theta_i f_i$

$$\|F^\theta\|_2 \leq \|F^\theta\|_1^{\frac{1}{4}} \|F^\theta\|_3^{\frac{3}{4}} \leq n^{\frac{1}{4}(\frac{1}{2}+\delta)}(Mn)^{\frac{3}{4}} = M^{\frac{3}{4}} n^{\frac{7}{8}+\frac{\delta}{4}}.$$

Thus, we can apply the Lemma for a set of functions $\{f_i/\|f_i\|_2\}_{i=1}^n$ in $L_2$ with $\beta = 3/8 + \delta/4$ and get

$$\|\sum_{i=1}^n a_i f_i\|_2 \leq C(M) n^{\frac{7+2\delta}{16}} \Big(\sum_{i=1}^n a_i^2\Big)^{\frac{1}{2}}.$$

Therefore, the function system $\{f_i/\|f_i\|_2\}_{i=1}^n$ satisfies the conditions **(a), (b)** (with another $M$) and $p = (7+2\delta)/16 < 1/2$. Let $\{\xi_i\}_1^n$ be the Rademacher functions and $m = 2^k \leq n$, now we can apply Theorem 1 (see (9) and the Remark to Theorem 1) to the random polynomial $F^\xi = \sum_1^n \xi_i f_i$ and get

$$\mathsf{P}\{\omega \in \Omega : \|F^\xi\|_{2^k,\infty} \leq C_1'(nk)^{\frac{1}{2}}\} \leq C_2' 2^{-qk}.$$

---

[2]Geometrically the Lemma implies that the convex hull of the set $B_\infty^d \cup (n^{1/2+\beta} \cdot B_1^d)$ has an inscribed sphere with radius of order $n^{1/4+\beta/2}$, here $B_\infty^d$ denotes the $d$-dimensional cube whose vertices have coordinates $\pm 1$, and $B_1^d = \{v \in \mathbb{R}^d : \sum_1^d |v_k| \leq 1\}$ (generalized octahedron).

Summing these inequalities up from $k = k_0$ to $\log n$ we get

$$\mathsf{P} \bigcup_{k=k_0}^{\log n} \{\omega \in \Omega : \|F^\xi\|_{2^k,\infty} \leq C_1'(nk)^{\frac{1}{2}}\} \leq C 2^{-qk_0}.$$

Thus, there exists a set of signs $\{\theta_i\}_{i=1}^n$ such that

$$\|F^\theta\|_{2^k,\infty} = \|\sum_{i=1}^n \theta_i f_i\|_{2^k,\infty} \geq C(M,\delta)\sqrt{nk} \tag{21}$$

for each $k = k_0, \ldots, \log n$. Moreover, we can take the constant $C(M,\delta)$ large enough in order to make the inequality (21) hold for all $k = 1, \ldots, \log n$.

Let us show now that (21) implies (18). Consider the set $e_k^* \subset X$, $\mu e_k^* = 2^{-k}$ such that

$$\int_{e_k^*} |F^\theta(x)| d\mu = \sup_{\substack{e \subset X \\ \mu e = 2^{-k}}} \int_e |\sum_{i=1}^n \theta_i f_i(x)| d\mu, \tag{22}$$

and let $\chi_{e_k^*}$ be its characteristic function. Using the triangle inequality for the integral-uniform norm and the fact that $\|g\|_{m,\infty} \leq m\|g\|_1$ for each $g$ (see (4)), obtain

$$\|F^\theta\|_{2^k,\infty} \leq \|F^\theta \cdot \chi_{e_k^*}\|_{2^k,\infty} + \|F^\theta \cdot (1 - \chi_{e_k^*})\|_{2^k,\infty} \leq$$
$$\leq 2^k \int_{e_k^*} |F^\theta(x)| d\mu + \sup_{x \notin e_k^*} |F^\theta(x)| \leq 2 \cdot 2^k \int_{e_k^*} |F^\theta(x)| d\mu.$$

This inequality along with (21) and (22) implies (18) for all $k = 1, \ldots, \log n$. To complete the proof of the Theorem it remains to prove the Lemma.

**Proof of the Lemma.** Consider a set of coefficients $\{a_i\}_1^n$ such that $\sum_1^n a_i^2 = n$ and $|a_i| < \sqrt{n}$. Define the following index sets

$$\sigma_k = \{j : 2^{-k}\sqrt{n} \leq |a_i| < 2^{-k+1}\sqrt{n}\} \qquad k \geq 1.$$

From Chebyshev's inequality it follows that $|\sigma_k| \leq 2^{2k}$. Define

$$w^{(k)} := \sum_{j \in \sigma_k} w_j.$$

From the triangle inequality it follows that $\|w^{(k)}\| \leq |\sigma_k| 2^{-k+1}\sqrt{n} \leq 2^{k+1}\sqrt{n}$. In order to estimate

$$W := \sum_{i=1}^n a_i w_i = \widetilde{w}^{(K)} + \sum_{k=1}^K w^{(k)}$$

notice that the residual term $\widetilde{w}^{(K)}$ belongs to the convex hull of the vectors $2^{-K}\sqrt{n} \sum_1^n \pm w_i$. So we get

$$\|\widetilde{w}^{(K)}\| \leq 2^{-K}\sqrt{n} \cdot C_{11} n^{\frac{1}{2}+\beta}.$$

Assembling all the facts we get

$$\|W\| \leq \sum_{k=1}^{K} \|w^{(k)}\| + \|\widetilde{w}^{(K)}\| \leq \sqrt{n} \sum_{k=1}^{K} 2^{k+1} + C_{11} 2^{-K} n^{1+\beta} < \sqrt{n}(2^{K+2} + C_{11} 2^{-K} n^{\frac{1}{2}+\beta}).$$

Choose $K \asymp (1/4 + \beta/2) \log n$ to obtain

$$\|W\| \leq C_{12}(C_{11}) n^{\frac{3}{4}+\frac{\beta}{2}},$$

which completes the proof of $(20')$ and Theorem 5.

To show that $(20')$ is precise consider a linear space with a basis $\{w_k\}_1^n$ and the norm

$$\|\sum_{k=1}^{n} a_k w_k\| := \max\Big(\sum_{k=1}^{[n^{1/2+\beta}]} |a_k|, \max_{k>n^{1/2+\beta}} |a_k|\Big).$$

Obviously, $\|w_k\| = 1$ and $\|\sum_1^n \pm w_k\| \leq n^{1/2+\beta}$. Consider a vector $W := \sum_{k=1}^{[n^{1/2+\beta}]} w_k$, for which $\|W\| = [n^{1/2+\beta}]$. Clearly, the inequality $(20')$ is precise for $W$.

## References

[1] R.N. Bhattacharya and R. Ranga Rao, *Normal Approximation and Asymptotic Expantions*, J. Wiley (1976)
[2] R. DeVore and G. Lorentz, *Constructive Approximation*, Springer Verlag, 1993
[3] P.G. Grigoriev, *Estimates for norms of random polynomials and their application*, Math. Notes **v. 69, N 6** (2001)
[4] J.-P. Kahane, *Some Random Series of Functions.* Heath mathematical monographs. Lexington, Mass., 1968.
[5] B. Kashin and L. Tzafriri, *Lower estimates for the supremum of some random processes*, East J. on Approx. **v. 1, N 1** (1995), 125–139.
[6] B. Kashin and L. Tzafriri, *Lower estimates for the supremum of some random processes, II*, East J. on Approx. **v. 1, N 3** (1995), 373–377.
[7] B. Kashin and L. Tzafriri, *Lower estimates for the supremum of some random processes, II*, Preprint, Max–Plank Institut für Mathematik, Bonn/95-85
[8] M. Ledoux and M. Talagrand, *Probability in Banach Spaces*, Springer Verlag, 1991.
[9] M. Marcus and G. Pisier, *Random Fourier Series with Applications to Harmonic Analysis*, Princeton Univ. Press, 1981.
[10] S. Montgomery-Smith and E.M. Semenov, *Embeddings of rearrangement invariant spaces that are not strictly singular,* Positivity **v. 4**, (2000), 397–402.
[11] S.M. Nikolskii, *Approximation of the Multivariate Functions and Embedding Theorems*, Nauka, Moscow, 1969 (in Russian)
[12] M. Riesz, *Formule d'interpolation pour la dérivée d'un polynome trigonométrique*, C. R. Acad. Sci. Paris 158 (1914), 1152–1154.
[13] V.I. Rotar', *Nonuniform estimate of the rate of convergence in multidimensional central limit theorem*, Theory of Probab. and Appl. **v. 15** (1970), 647–665
[14] R. Salem, A. Zygmund, *Some properties of trigonometric series whose terms have random signs*, Acta Math. **91**(1954) 245–301.